\documentclass[article,12pt]{amsart}
\usepackage{bbm}
\usepackage{amsmath}
\usepackage{txfonts}
\usepackage{stmaryrd}
\usepackage{amssymb}
\usepackage{amsfonts}
\usepackage{times}

\newcommand{\D}{{\mathbb D}}

\def\D(a,r){\mathbb D(a,r)}
\def\D\m\D(a,r_0){\Delta\setminus\Delta(a,r_0)}

\def\D{\mathbb D}

\newcommand{\DD}{{\mathbb D}}
\newcommand{\CC}{{\mathbb C}}

\def\a{\alpha}

\def\msk{\medskip}

   \def\a{\alpha} 
\def\b{\beta}   
 
\def\qand{\quad\mbox{ and }\quad}

\def\bege{\begin{equation}} \def\ende{\end{equation}}
\def\begr{\begin{eqnarray}} \def\endr{\end{eqnarray}}
\def\BB{ \mathcal{B} }

\numberwithin{equation}{section}

\begin{document}

\title[Essential norm  of    weighted composition operators]{Essential norm and a new characterization of    weighted composition operators from weighted Bergman spaces and Hardy spaces into the Bloch space}

\author{Songxiao Li, Ruishen Qian  and Jizhen Zhou}

{\footnotesize
\address{Songxiao Li\\ Institute of Systems Engineering, Macau University of Science and Technology, Avenida Wai Long, Taipa, Macau. }    \email{jyulsx@163.com}
\address{Ruishen Qian\\ School of Mathematics and Computation Science,
 Lingnan Normal University, Zhanjiang 524048, Guangdong, P. R. China}          \email{qianruishen@sina.cn}
\address{Jizhen Zhou\\ School of Sciences, Anhui University of Science and Technology, Huainan, Anhui 232001, P. R. China}   \email{hope189@163.com}}


\begin{abstract}
In this paper, we give some estimates for the essential norm and a new characterization for the boundedness and compactness  of  weighted composition operators from weighted Bergman spaces and Hardy spaces to the Bloch  space.
\end{abstract}
\maketitle

\section{Introduction}

Let $\mathbb{D}$ be the open unit disk in the
complex plane $\CC$ and $H(\mathbb{D})$ be the space of analytic
functions on $\DD$. For $0<p<\infty$ and $\alpha>-1$, the weighted Bergman space,
denoted by $A^p_\alpha$,  is the set of all functions $f \in H(\DD)$
satisfying
$$
\|f\|_{A_\alpha^p}^p=(\alpha+1)\int_{\DD} |f(z)|^p (1-|z|^2)^\alpha
dA(z) <\infty,
$$
where $dA$ is the normalized Lebesgue area measure in $\DD$ such that
$A(\DD)= 1$. The Hardy space $H^p$ is the space consisting of all  $f\in
H(\DD)$ such that
$$\|f\|_{H^p}^p=\sup_{0<r<1}\frac{1}{2\pi}\int_0^{2\pi}|f(re^{i\theta})|^pd\theta<\infty.$$

 The Bloch space, denoted by $\mathcal{B}=\mathcal{B}(\DD)$,  is the space of all $f \in
H(\DD)$ such that
  \begr \|f\|_\b  = \sup_{z \in \DD}(1-|z|^2) |f'(z)|<\infty. \nonumber\endr Under the norm
$\|f\|_{\mathcal{B}}=|f(0)|+ \|f\|_\b$, the Bloch space is a Banach
space. See \cite{zhu} for more information of the Bloch space.

Let $v:\DD\rightarrow R_+$ be a continuous, strictly positive and bounded function. An $f\in H(\DD)$ is said to belong to the weighted space, denoted by
$H^\infty_v$, if
$$
\|f\|_v=\sup_{z \in \mathbb{D}}v(z)|f(z)|<\infty.
$$
 $H^\infty_v$ is a Banach space with the norm $\| \cdot\|_v$. The weighted $v$ is called radial, if $v(z)=v(|z|)$ for all $z\in \DD$. For a weight $v$, the associated weight  $\widetilde{v}$  is defined as follows.
$$
\widetilde{v}=(\sup\{|f(z)|:  f\in H^\infty_v, \|f\|_v \leq 1  \})^{-1},   z\in \DD.
$$
When $v=v_\a(z)=(1-|z|^2)^\alpha(0<\a<\infty)$, it is easy to check that $\widetilde{v}_\a(z)=v_\a(z)$. In this case, we denote $H^\infty_v$ by $H^\infty_{v_\a}$
 and $\|f\|_{v_\a}=\sup_{z \in \mathbb{D}} |f(z)|(1-|z|^2)^\alpha $.

Let $S(\DD)$ denote the set of all analytic self-maps  of $\mathbb{D}$.  Let $u \in H(\mathbb{D})$ and $\varphi\in S(\DD)$. For $f \in H(\mathbb{D})$, the   composition operator $C_\varphi$ and the  multiplication operator $M_u$ are defined by
$$
(C_\varphi f)(z) =    f(\varphi(z)) ~~~\mbox{and}~~~~~(M_u f)(z)=u(z)f(z),$$
respectively. The weighted composition operator $uC_\varphi$   is defined by
$$
(uC_\varphi f)(z) =u(z)\cdot f(\varphi(z)), \  \  f \in H(\mathbb{D}). $$
It is clear that the weighted composition operator $uC_\varphi$ is the generalization of
$C_\varphi$ and $M_u$.  A basic and interesting problem concerning  concrete operators (such as  composition operator, multiplication operator, Volterra operator, Toeplitz operator, Hankel operator and other integral type operators) is to relate operator theoretic properties to their function theoretic properties of their symbols, which attracted a lot of attention recently,
the reader can refer to \cite{cm} and \cite{zhu}.

  It is  well known   that $C_\varphi$ is bounded on $\mathcal{B}$ by the Schwarz-Pick
lemma for any $\varphi\in S(\DD)$. The compactness of $C_\varphi$   on
$\mathcal{B}$ was studied in for example \cite{mm, t, wzz}.  In \cite{wzz},
Wulan, Zheng and Zhu proved that, for any $\varphi \in S(\DD)$,  $C_\varphi  :\mathcal{B}\rightarrow
\mathcal{B}$ is compact if and only if $\lim_{j\rightarrow\infty}\| \varphi^j \|_\mathcal{B}=0.$
This result has been generalized to Bloch-type spaces by  Zhao in
\cite{zhao} and shows that $C_\varphi
 :\mathcal{B}^\alpha\rightarrow \mathcal{B}^\beta$ is compact if and only if
$\lim_{j\rightarrow\infty} j^{\alpha-1} \|\varphi^j
\|_{\mathcal{B}^\beta}=0.$   For some results on composition operator amd related operators mapping into the Bloch space see, for example, \cite{c, col,ls888, ls00, ls1, lz, lou, mz1, mm, mz2, osz,   st5, st6, x,   y2, yz, zhao, zxl} and the related references therein.

In \cite{ls888}, the first author of this paper and Stevi\'c obtained a characterization of the boundedness and compactness of weighted composition operator $uC_\varphi:A^p_\a \to \mathcal{B} $. Among others, we proved the following result.\msk

\noindent{\bf Theorem A. } {\it Let  $1\leq p<\infty$, $\a>-1$, $u \in H(\mathbb{D})$ and $\varphi\in S(\DD)$ such that $uC_\varphi:A^p_\a \to \mathcal{B} $ is
bounded.    Then $uC_\varphi:A^p_\a \to \mathcal{B} $ is compact if and only if
 $$\lim_{|\varphi(z)|\rightarrow
1}\frac{(1-|z|^{2})|u'(z)|}{(1-|\varphi(z)|^{2})^{(2+\a)/p}}=0\qand
\lim_{|\varphi(z)|\rightarrow
1}\frac{(1-|z|^{2})|u(z)\varphi'(z)|}{(1-|\varphi(z)|^{2})^{ (2+\a+p)/p}}=0.$$    }

In \cite{col}, Colonna obtained a new characterization by using two families functions, among others, she obtained the following result.\msk

\noindent{\bf Theorem B. }  {\it Let  $1\leq p<\infty$, $\a>-1$, $u \in H(\mathbb{D})$ and $\varphi\in S(\DD)$ such that $uC_\varphi:A^p_\a \to \mathcal{B} $ is
bounded.    Then $uC_\varphi:A^p_\a \to \mathcal{B} $ is compact if and only if
  $$\lim_{|a|\to 1}\big\|uC_\varphi  f_a   \big\|_{\mathcal{B}  }=0
~~~~~ \mbox{and} ~~~~~ \lim_{|a|\to 1}\big\|uC_\varphi  g_a \big\|_{\mathcal{B}  }=0\nonumber , $$
where
$$
f_a(z)=\frac{(1-|a|^2)^{1+(2+\a)(1-1/p)}}{(1-\overline{a}z)^{3+\alpha}},  g_a(z)=\frac{(1-|a|^2)^{1+(2+\a)(1-1/p)+1/p}}{(1-\overline{a}z)^{3+\alpha+1/p}}.
$$ }

 In \cite{col}, Colonna also obtained two characterizations for the   compactness of weighted composition operator $uC_\varphi:H^p \to \mathcal{B} $.\msk

\noindent{\bf Theorem C. }  {\it Let  $1\leq p<\infty$,  $u \in H(\mathbb{D})$ and $\varphi\in S(\DD)$ such that $uC_\varphi:H^p \to \mathcal{B} $ is
bounded.    Then the following statements are equivalent:

(a) $uC_\varphi:H^p \to \mathcal{B} $ is compact.

(b)  $$\lim_{|a|\to 1}\big\|uC_\varphi  p_a
   \big\|_{\mathcal{B}  }=0
~~~~~ \mbox{and} ~~~~~ \lim_{|a|\to 1}\big\|uC_\varphi  q_a \big\|_{\mathcal{B}  }=0\nonumber , $$
where
$$
p_a(z)=\frac{(1-|a|^2)^{2-1/p}}{(1-\overline{a}z)^2},  q_a(z)=\frac{(1-|a|^2)^2}{(1-\overline{a}z)^{2+1/p}}.
$$

(c) $$\lim_{|\varphi(z)|\rightarrow
1}\frac{(1-|z|^{2})|u'(z)|}{(1-|\varphi(z)|^{2})^{1/p}}=0\qand
\lim_{|\varphi(z)|\rightarrow
1}\frac{(1-|z|^{2})|u(z)\varphi'(z)|}{(1-|\varphi(z)|^{2})^{ (1+p)/p}}=0.$$   }

The purpose of this paper is to give some estimates for the essential norm of  the operator $uC_{\varphi}:A^p_\alpha \rightarrow \BB $(as well as $uC_{\varphi}:H^p \rightarrow \BB $), in particular, by using $\|uC_\varphi  f_a  \|_{\mathcal{B}  }$ and $\|uC_\varphi  g_a  \|_{\mathcal{B}  }$(as well as $\|uC_\varphi  p_a  \|_{\mathcal{B}  }$ and $\|uC_\varphi  q_a  \|_{\mathcal{B}  }$). Moreover, we give a new characterization for the boundedness, compactness and essential norm of  the operator $uC_{\varphi}:A^p_\alpha \rightarrow \BB $(as well as $uC_{\varphi}:H^p \rightarrow \BB $) by using $\varphi^j$.

Recall that the essential norm of a bounded linear operator $T:X\rightarrow Y$ is its distance to the set of
compact operators $K$ mapping $X$ into $Y$, that is,
$$\|T\|_{e, X\rightarrow Y}=\inf\{\|T-K\|_{X\rightarrow Y}: K~\mbox{
is compact}~~\},$$ where  $X,Y$ are Banach spaces and
$\|\cdot\|_{X\rightarrow Y}$ is the operator norm.

Throughout this paper,  we say that $A\lesssim B$ if there exists a constant
$C$ such that $A\leq CB$. The symbol $A\approx B$ means that
$A\lesssim B\lesssim A$.

\section{ Essential norm of $uC_\varphi $}

In this section, we give two estimates for the essential norm of the operator $uC_\varphi:A^p_\a \to \mathcal{B} $ and the operator $uC_\varphi:H^p\to \mathcal{B} $, respectively.  \msk

\noindent {\bf Theorem 2.1.} {\it Let  $1\leq p<\infty$, $\a>-1$, $u \in H(\mathbb{D})$ and $\varphi\in S(\DD)$ such that $uC_\varphi:A^p_\a \to \mathcal{B} $ is
bounded.    Then
  \begr  \|uC_{\varphi}\|_{e,A^p_\a \to \BB } \approx\max \Big\{ A,  B \Big\}\approx \max \Big\{ P,  Q \Big\}   , \nonumber
\endr
where
$$
A:=\limsup_{|a|\to
1 }\left\|uC_{\varphi}\left(f_a\right)\right\|_{\BB },   ~~~~~~~~~~B:= \limsup_{|a|\to
1 }\left\|uC_{\varphi}\left(g_a \right)\right\|_{\BB },
$$
$$
P:=\limsup_{|\varphi(z)|\rightarrow1}\frac{(1-|z|^{2})|u'(z)|}{(1-|\varphi(z)|^{2})^{(2+\a)/p}}
, ~~~~  Q:= \limsup_{|\varphi(z)|\rightarrow1}\frac{(1-|z|^{2})|u(z)\varphi'(z)|}{(1-|\varphi(z)|^{2})^{ (2+\a+p)/p}}  .
$$}

\begin{proof}  First we prove that
$$ \max \Big\{ A,   B \Big\} \lesssim  \|uC_\varphi\|_{e,A^p_\a\to \BB } .
$$
 Let $a\in \DD$. Define
 $$
f_a(z)=\frac{(1-|a|^2)^{1+(2+\a)(1-1/p)}}{(1-\overline{a}z)^{3+\alpha}},~~~~~~ g_a(z)=\frac{(1-|a|^2)^{1+(2+\a)(1-1/p)+1/p}}{(1-\overline{a}z)^{3+\alpha+1/p}} , ~~z\in \DD.\nonumber
$$
It is easy to check that $f_a, g_a\in A^p_\a$ and $\|f_a\|_{A^p_\a}  \lesssim 1,   \|g_a\|_{A^p_\a} \lesssim 1$ for all $a\in \DD$ and
$f_a, g_a$ converges   to zero uniformly on compact subsets of $\DD$ as $|a|\to
1$.   Thus, for any compact operator $K:A^p_\a\to \BB$, by Lemma 3.7 of \cite{t2}  we have
 $$
 \lim_{|a|\to 1}\|Kf_a\|_{\BB }=0, ~~~\lim_{|a|\to 1}\|Kg_a\|_{\BB }=0.
 $$
  Hence
   \begr
 \|uC_\varphi-K\|_{A^p_\a\to \BB}&\gtrsim &\|(uC_\varphi-K)f_a\|_{\BB } \geq   \|uC_\varphi f_a \|_{\BB }-\|Kf_a\|_{\BB }, \nonumber
\endr
and
  \begr
 \|uC_\varphi-K\|_{A^p_\a\to \BB}&\gtrsim &\|(uC_\varphi-K)g_a\|_{\BB } \geq   \|uC_\varphi g_a \|_{\BB}-\|Kg_a\|_{\BB}.\nonumber
\endr
Taking $\limsup_{|a|\to 1}$ to the last two inequalities on both sides,  we obtain
   \begr  \|uC_\varphi-K\|_{ A^p_\a\to \BB}  \gtrsim  A, ~~~ ~~~~~\|uC_\varphi-K\|_{A^p_\a\to \BB} \gtrsim B .\nonumber
\endr
Therefore, by the definition of the essential norm, we get
   \begr   \|uC_\varphi\|_{e,A^p_\a\to \BB }= \inf_{K}  \|uC_\varphi-K\|_{A^p_\a\to \BB }  \gtrsim \max \Big\{ A,   B \Big\} .\nonumber
\endr
Next, set
$$
h_a(z)=f_a-g_a, ~~~~~~~~~~~~~~~~ k_a(z)=f_a-\frac{3+\a}{3+\a+1/p}g_a.
$$
It is also easy to check that $h_a, k_a\in A^p_\a$ and $\|h_a\|_{A^p_\a} \lesssim  1,   \|k_a\|_{A^p_\a} \lesssim 1$ for all $a\in \DD$ and
$h_a, k_a$  converges to zero uniformly on compact subsets of $\DD$ as $|a|\to
1$.     Hence, for any $b_j\in \DD$ such that $|\varphi(b_j)|\rightarrow 1$ and  any compact operator $K:A^p_\a\to \BB$, we have
   \begr
 \|uC_\varphi-K\|_{A^p_\a\to \BB}&\gtrsim &\|(uC_\varphi-K)h_{\varphi(b_j)}\|_{\BB } \geq   \|uC_\varphi h_{\varphi(b_j)} \|_{\BB }-\|Kh_{\varphi(b_j)}\|_{\BB }, \nonumber
\endr
and
  \begr
 \|uC_\varphi-K\|_{A^p_\a\to \BB}&\gtrsim &\|(uC_\varphi-K)k_{\varphi(b_j)}\|_{\BB } \geq   \|uC_\varphi k_{\varphi(b_j)} \|_{\BB}-\|K k_{\varphi(b_j)}\|_{\BB}.\nonumber
\endr
Taking $\limsup_{|{\varphi(b_j)}|\to 1}$ to the last two inequalities on both sides   we get
 \begr
 \|uC_\varphi-K\|_{A^p_\a\to \BB} \gtrsim  \limsup_{|\varphi(b_j)|\rightarrow1} \|uC_\varphi h_{\varphi(b_j)} \|_{\BB }  \gtrsim  \limsup_{|\varphi(b_j)|\rightarrow1} \frac{(1-|b_j|^{2})|u'(b_j)|}{(1-|\varphi(b_j)|^{2})^{(2+\a)/p}} =P  , \nonumber
\endr
and
  \begr
 \|uC_\varphi-K\|_{A^p_\a\to \BB} \gtrsim  \limsup_{|\varphi(b_j)|\rightarrow1}  \|uC_\varphi k_{\varphi(b_j)} \|_{\BB} \gtrsim  \limsup_{|\varphi(b_j)|\rightarrow1}   \frac{(1-|b_j|^{2})|u(b_j)\varphi'(b_j)|}{(1-|\varphi(b_j)|^{2})^{ (2+\a+p)/p}}=Q.\nonumber
\endr
By the definition of the essential norm, we obtain
   \begr   \|uC_\varphi\|_{e,A^p_\a\to \BB }= \inf_{K}  \|uC_\varphi-K\|_{A^p_\a\to \BB }\gtrsim \max \Big\{ P,   Q \Big\} .\nonumber
\endr

Finally, we prove that $$ \|uC_\varphi\|_{e,A^p_\a\to \BB }  \lesssim \max \Big\{ A,   B \Big\},~~~~ \mbox{and}~~~~ \|uC_\varphi\|_{e,A^p_\a\to \BB }  \lesssim \max \Big\{ P,   Q \Big\}. $$
For $r\in [0,1)$,  set $K_r: H(\D)\to H(\D)$ by $$(K_r f)(z)=f_r(z)=f(rz), ~~f\in H(\D).$$
 It is clear that  $K_r$ is compact on
$A^p_\a$ and $ \|K_r\|_{A^p_\a\to A^p_\a}\leq 1.$   Let $\{r_j\}\subset (0,1)$ be a sequence such that
$r_j\to 1$ as $j\to\infty$. Then for all positive integer
$j$, the operator $uC_\varphi K_{r_j}: A^p_\a\rightarrow \BB $ is   compact. By the definition of the essential norm we have
  \begr
\|uC_\varphi\|_{e,A^p_\a\to \BB } \leq \limsup_{j\to\infty}\|uC_\varphi- uC_\varphi K_{r_j}\|_{A^p_\a\to
\BB }.
\endr
  Thus, we only need to show that
    \begr && \limsup_{j\to\infty}\|uC_\varphi-uC_\varphi K_{r_j}\|_{ A^p_\a\to \BB }   \lesssim \max \Big\{ A,   B \Big\},
\endr
and  \begr && \limsup_{j\to\infty}\|uC_\varphi-uC_\varphi K_{r_j}\|_{ A^p_\a\to \BB }   \lesssim \max \Big\{ P,   Q \Big\}.
\endr
 For any $f\in A^p_\a$ such that $\|f\|_{A^p_\a}\leq 1$, we consider
 \begr &&\|( uC_\varphi- uC_\varphi K_{r_j})f\|_{\BB }=|u(0)f(\varphi(0))-u(0)f(r_j\varphi(0))|+\|u\cdot (f-f_{r_j})\circ\varphi\|_\b.  \nonumber
\endr
It is clear that
 $\lim_{j\to\infty}|u(0)f(\varphi(0))-u(0)f(r_j\varphi(0))|=0.$
  Now we estimate
  \begr
&& \limsup_{j\to\infty}\|u\cdot (f-f_{r_j}) \circ\varphi\|_\beta\nonumber\\
&=& \limsup_{j\to\infty}\sup_{|\varphi(z)|\leq r_N}(1-|z|^2)|(f-f_{r_j})'(\varphi(z))||\varphi'(z)||u(z)|\nonumber\\
&  & + \limsup_{j\to\infty}\sup_{|\varphi(z)|> r_N}(1-|z|^2) |(f-f_{r_j})'(\varphi(z))||\varphi'(z)||u(z)| \nonumber\\
& & +\limsup_{j\to\infty}\sup_{|\varphi(z)|\leq r_N}(1-|z|^2) |(f-f_{r_j})(\varphi(z))||u'(z)| \nonumber\\
&  & + \limsup_{j\to\infty}\sup_{|\varphi(z)|> r_N}(1-|z|^2) |(f-f_{r_j})(\varphi(z))||u'(z)|\nonumber\\
&=&Q_1+Q_2+Q_3+Q_4,
\endr
 where $N\in \mathbb{N }$ is large enough such that $r_j\geq\frac{1}{2}$ for all $j\geq N$,
$$
Q_1:=\limsup_{j\to\infty}\sup_{|\varphi(z)|\leq r_N}(1-|z|^2)|(f-f_{r_j})'(\varphi(z))||\varphi'(z)||u(z)| ,
$$
$$
Q_2:=\limsup_{j\to\infty}\sup_{|\varphi(z)|> r_N}(1-|z|^2)|(f-f_{r_j})'(\varphi(z))||\varphi'(z)||u(z)|,
$$
$$
Q_3:=\limsup_{j\to\infty}\sup_{|\varphi(z)|\leq r_N}(1-|z|^2)|(f-f_{r_j})(\varphi(z))||u'(z)|,
$$
and
$$
Q_4:=\limsup_{j\to\infty}\sup_{|\varphi(z)|> r_N}(1-|z|^2)|(f-f_{r_j})(\varphi(z))||u'(z)|.
$$
 Since $uC_\varphi:A^p_\a\to \mathcal{B} $ is bounded, applying the operator $uC_\varphi$ to 1 and $z$, we easily get  that $u\in  \mathcal{B} $ and
$$
\widetilde{K}:=\sup_{z\in \DD}(1-|z|^2) |\varphi'(z)||u(z)|<\infty.
$$   Since   $r_jf'_{r_j}\to f'$ uniformly on compact subsets of $\DD$ as $j\to\infty$, we have
  \begr
Q_1  \leq \widetilde{ K} \limsup_{j\to\infty}\sup_{|w|\leq r_N}|f'(w)- r_jf' (r_j w)| =0.
\endr
Also, from the fact that $u \in  \mathcal{B} $ and $f_{r_j}\to f$ uniformly on compact subsets of $\DD$ as $j\to\infty$, we have
  \begr
Q_3  \leq  \|u\|_{\BB }\limsup_{j\to\infty}\sup_{|w|\leq r_N}|f(w)- f(r_j w)| =0.
\endr
Next we consider $Q_2$.  We have $Q_2\leq\limsup_{j\to\infty}(S_1+S_2), $  where
 \begr
S_1:=\sup_{|\varphi(z)|>r_N}(1-|z|^2)  |f'(\varphi(z))||\varphi'(z)||u(z)| \nonumber \endr
  and
 \begr
S_2:=\sup_{|\varphi(z)|>r_N}(1-|z|^2)  r_j|f'(r_j\varphi(z))||\varphi'(z)||u(z)|.\nonumber
\endr
 First we estimate $S_1$. Using the fact that $\|f\|_{A^p_\a}\leq 1$, we have
    \begr S_1&=&\sup_{|\varphi(z) |>r_N}(1-|z|^2)|f'(\varphi(z))| |\varphi'(z)||u(z)|  \nonumber\\
&\lesssim &\frac{1}{ r_N }   \|f\|_{A^p_\a}\sup_{|\varphi(z) |>r_N}(1-|z|^2)   |\varphi'(z)||u(z)|  \frac{ |\varphi(z)|  }{(1-|\varphi(z)|^2)^{ (2+\a+p)/p}}\nonumber\\
&\lesssim & \frac{1}{p}  \sup_{|\varphi(z) |>r_N}\sup_{|a|>r_N}(1-|z|^2)    |\varphi'(z)||u(z)|  \times \frac{|\varphi(z)| }{(1-|\varphi(z)|^2)^{ (2+\a+p)/p}} \nonumber\\
&\lesssim&   \sup_{|a|>r_N} \left\|uC_\varphi (f_a- g_a)\right\|_{\BB}  \nonumber\\
& \lesssim&   \sup_{|a|>r_N} \left\|uC_\varphi f_a\right\|_{\BB }+   \sup_{|a
|>r_N} \left\|uC_\varphi g_a\right\|_{\BB } .\nonumber
\endr
 Taking limit
as $N\to\infty$ we obtain \begr
\limsup_{j\to\infty}S_1&\lesssim&  \limsup_{|a|\to 1} \frac{ (1-|z|^2) |\varphi'(z)||u(z)|   }{(1-|\varphi(z)|^2)^{(2+\a+p)/p}}=Q\nonumber\\
&\lesssim & \limsup_{|a|\to 1}\left \|uC_\varphi f_a\right\|_{\BB }+   \limsup_{|a|\to1}\left\|uC_\varphi g_a\right\|_{\BB }.\nonumber
\endr
Similarly, we have
  \begr
\limsup_{j\to\infty}S_2& \lesssim &  \limsup_{|a|\to 1} \frac{(1-|z|^2)  |\varphi'(z)||u(z)|   }{(1-|\varphi(z)|^2)^{ (2+\a+p)/p}}=Q\nonumber\\
&\lesssim & \limsup_{|a|\to 1}\left \|uC_\varphi f_a\right\|_{\BB }+   \limsup_{|a|\to 1}\left\|uC_\varphi g_a\right\|_{\BB }.\nonumber
\endr
{\it i.e.}, we get
that \begr
  Q_2\lesssim  Q\lesssim   A+   B \lesssim \max \Big\{ A,   B \Big\}  .
\endr

Next we consider $Q_4$.  We have $Q_4\leq\limsup_{j\to\infty}(S_3+S_4), $  where
 \begr
S_3:=\sup_{|\varphi(z)|>r_N}(1-|z|^2) |f(\varphi(z))||u'(z)|,
S_4:=\sup_{|\varphi(z)|>r_N}(1-|z|^2) |f(r_j\varphi(z))||u'(z)|.\nonumber
\endr
 Similarly,   we have
    \begr S_3&\lesssim &  \sup_{|\varphi(z) |>r_N}\sup_{|a|>r_N}(1-|z|^2)  |u'(z)|  \frac{ 1}{(1-|\varphi(z)|^2)^{(2+\a)/p}} \nonumber\\
&\lesssim &   \sup_{|a|>r_N} \left\|uC_\varphi f_a- \frac{3+\a}{3+\a+1/p} uC_\varphi g_a \right\|_{\BB}  \nonumber\\
&\leq&  \sup_{|a|>r_N} \left\|uC_\varphi f_a\right\|_{\BB }+ \frac{3+\a}{3+\a+1/p}\sup_{|a|>r_N} \left\|uC_\varphi g_a\right\|_{\BB } \nonumber\\
&\leq &\sup_{|a|>r_N} \left\|uC_\varphi f_a\right\|_{\BB }+ \sup_{|a|>r_N} \left\|uC_\varphi g_a\right\|_{\BB } .\nonumber
\endr
 Taking limit
as $N\to\infty$ we obtain
 \begr
\limsup_{j\to\infty}S_3&\lesssim & \limsup_{|a|\to 1 }  \frac{ (1-|z|^2)  |u'(z)|  }{(1-|\varphi(z)|^2)^{(2+\a)/p}} =P\nonumber\\
&\lesssim &\limsup_{|a|\to 1 }\left\|uC_\varphi f_a\right\|_{\BB }+ \limsup_{|a|\to
1 }\left\|uC_\varphi g_a\right\|_{\BB }= A+B.  \nonumber
\endr
Similarly, we have  $\limsup_{j\to\infty}S_4  \lesssim P \lesssim A+ B,$  {\it i.e.}, we get
that
 \begr   Q_4 \lesssim P \lesssim A  +  B.
\endr
  Hence, by (2.4), (2.5), (2.6), (2.7) and (2.8) we get
   \begr
 \limsup_{j\to\infty}\|uC_\varphi- uC_\varphi K_{r_j}\|_{ A^p_\a\to \BB }
 &=&\limsup_{j\to\infty}\sup_{ \|f\|_{A^p_\a} \leq 1}\|(uC_\varphi- uC_\varphi K_{r_j})f\|_{\BB } \nonumber\\
&=&\limsup_{j\to\infty}\sup_{ \|f\|_{A^p_\a} \leq 1} \|u\cdot (f-f_{r_j}) \circ\varphi\|_{\beta }\nonumber\\
  &\lesssim & P+Q \lesssim  A  +  B  .
\endr
 Therefore, by (2.1)  and (2.9), we obtain
  \begr \|uC_\varphi\|_{e,A^p_\a\to \BB }  \lesssim P+Q \lesssim\max \Big\{ P,   Q \Big\}     \nonumber
\endr
and  \begr \|uC_\varphi\|_{e,A^p_\a\to \BB }  \lesssim A+B \lesssim   \max \Big\{ A,   B \Big\} . \nonumber
\endr
This completes the proof of the theorem. \end{proof} \msk

The Hardy space $H^p$ can be viewed as limiting space of $A^p_\a$ as $\a$ decreases to $-1$. Hence, from Theorem 2.1, we get the following result.\msk

\noindent {\bf Theorem 2.2.} {\it Let  $1\leq p<\infty$,  $u \in H(\mathbb{D})$ and $\varphi\in S(\DD)$ such that $uC_\varphi:H^p\to \mathcal{B} $ is
bounded.    Then
  \begr  \|uC_{\varphi}\|_{e,H^p \to \BB }& \approx& \max \Big\{ \limsup_{|a|\to
1 }\left\|uC_{\varphi}(p_a)\right\|_{\BB },    \limsup_{|a|\to
1 }\left\|uC_{\varphi}(q_a )\right\|_{\BB } \Big\}\nonumber\\
&\approx& \max \Big\{\limsup_{|\varphi(z)|\rightarrow1}\frac{(1-|z|^{2})|u'(z)|}{(1-|\varphi(z)|^{2})^{1/p}}  , \limsup_{|\varphi(z)|\rightarrow1}\frac{(1-|z|^{2})|u(z)\varphi'(z)|}{(1-|\varphi(z)|^{2})^{ (1+p)/p}}  \Big\} . \nonumber
\endr  }

From Theorems 2.1 and 2.2, we immediately get the following two corollaries.\msk

\noindent {\bf Corollary 2.1.} {\it Let  $1\leq p<\infty$, $\a>-1$  and $\varphi\in S(\DD)$ such that $C_\varphi:A^p_\a \to \mathcal{B} $ is
bounded.    Then
  \begr  \|C_{\varphi}\|_{e,A^p_\a \to \BB } &\approx&  \limsup_{|a|\to
1 }\left\|C_{\varphi}\left(f_a\right)\right\|_{\BB } \approx \limsup_{|a|\to
1 }\left\|C_{\varphi}\left(g_a \right)\right\|_{\BB } \nonumber\\
 &\approx & \limsup_{|\varphi(z)|\rightarrow1}\frac{(1-|z|^{2})|\varphi'(z)|}{(1-|\varphi(z)|^{2})^{ (2+\a+p)/p}}  . \nonumber
\endr
 }

\noindent {\bf Corollary 2.2.} {\it Let  $1\leq p<\infty$  and $\varphi\in S(\DD)$  such that $C_\varphi:H^p \to \mathcal{B} $ is
bounded.    Then
  \begr  \|C_{\varphi}\|_{e,H^p \to \BB } &\approx&  \limsup_{|a|\to
1 }\left\|C_{\varphi}\left(p_a\right)\right\|_{\BB } \approx \limsup_{|a|\to
1 }\left\|C_{\varphi}\left(q_a \right)\right\|_{\BB } \nonumber\\
 &\approx & \limsup_{|\varphi(z)|\rightarrow1}\frac{(1-|z|^{2})|\varphi'(z)|}{(1-|\varphi(z)|^{2})^{ (1+p)/p}}  . \nonumber
\endr
 }

\section{New characterization of $uC_\varphi $}

In this section, motivated by \cite{el}, we give a new characterization for the boundedness, compactness and essential norm for the weighted composition operators  $uC_\varphi: A^p_\alpha  \to \mathcal{B}  $ and $uC_\varphi: H^p \to \mathcal{B}  $. For this purpose,  we state some lemmas which will be used. \msk

\noindent{\bf Lemma 3.1.} \cite{mon} {\it Let $v$ and $w$ be radial, non-increasing weights tending to zero at the boundary of $\DD$.  Then the following statements hold.

(a) The weighted composition operator $uC_\varphi:H_{v}^{\infty}\rightarrow H_{w}^{\infty}$ is bounded if and only if
$$
\sup_{z\in \DD}\frac{w(z)}{\widetilde{v}(\varphi(z))}|\varphi(z)|<\infty.
$$
Moreover, the following holds
  \begr\|uC_\varphi\|_{H_{v}^{\infty}\rightarrow H_{w}^{\infty}}=\sup_{z\in \DD}\frac{w(z)}{\widetilde{v}(\varphi(z))}|\varphi(z)|.\nonumber
\endr

(b) Suppose $uC_\varphi:H_{v}^{\infty}\rightarrow H_{w}^{\infty}$ is bounded. Then
  \begr\|uC_\varphi\|_{e, H_{v}^{\infty}\rightarrow H_{w}^{\infty}}=\lim_{s\to 1^-}\sup_{|\varphi(z)|>s}\frac{w(z)}{\widetilde{v}(\varphi(z))}|\varphi(z)|.\nonumber
\endr}

\noindent{\bf Lemma 3.2.} \cite{h2} {\it Let $v$ and $w$ be radial, non-increasing weights tending to zero at the boundary of $\DD$. Then the following statements hold.

(a)    $uC_\varphi:H_{v}^{\infty}\rightarrow H_{w}^{\infty}$ is bounded if and only if
$$\sup_{k\geq 0}\frac{\|u \varphi^k\|_w}{\|z^k\|_v}<\infty,$$
with the norm comparable to the above supermum.

(b)  Suppose $uC_\varphi:H_{v}^{\infty}\rightarrow H_{w}^{\infty}$ is bounded. Then
  \begr\|uC_\varphi\|_{e,H_{v}^{\infty}\rightarrow H_{w}^{\infty} }=\limsup_{k\to
\infty}\frac{\|u \varphi^k\|_w}{\|z^k\|_v}.\nonumber
\endr}

 \noindent{\bf Lemma 3.3.}\cite{h1}  {\it For $\a>0$, we have $ \lim_{k\rightarrow\infty}k^\a\|z^{k-1}\|_{v_\a}=(\frac{2\a}{e})^\a$. }\msk

\noindent{\bf Theorem 3.1.} {\it   Let    $1 \leq p <\infty$, $\a>-1$, $u \in H(\mathbb{D})$ and   $\varphi\in S(\DD)$. Then the operator $uC_\varphi: A^p_\alpha  \to \mathcal{B}   $ is bounded if and only if
  \begr\sup_{j\geq 1} j^{(2+\a)/p} \|I_u(\varphi^j)\|_{\mathcal{B} } <\infty~~~~~~~\mbox{and}~~~~~~
 \sup_{j\geq 1} j^{(2+\a)/p}\|J_u(\varphi^{j-1})\|_{\mathcal{B} }   <\infty, \endr
where
$$
I_ug(z)=\int_0^z g'(\xi)u(\xi)d\xi, ~~~~~~~~~~~~J_ug(z)=\int_0^z g(\xi)u'(\xi)d\xi,   ~~z\in \DD,  g\in H(\DD).
$$ }

  \proof
  By Theorem A, $uC_\varphi: A^p_\alpha  \to \mathcal{B}   $ is bounded if and only if
\begr ~~~~\sup_{z\in \DD}\frac{(1-|z|^{2})|u'(z)|}{(1-|\varphi(z)|^{2})^{(2+\a)/p}}< \infty~~\mbox{and}~~
\sup_{z\in \DD}\frac{(1-|z|^{2})|u(z)\varphi'(z)|}{(1-|\varphi(z)|^{2})^{ (2+\a+p)/p}}<\infty,
\endr
 which are equivalent to the weighted composition operator $u'C_\varphi: H^\infty_{v_{(2+\a)/p}}\rightarrow  H^\infty_{v_1}$ is bounded and $u\varphi'C_\varphi: H^\infty_{v_{(2+\a+p)/p}}\rightarrow  H^\infty_{v_1}$ is bounded, respectively.
By Lemma 3.2, we see that two inequalities in (3.2) are equivalent to
$$
\sup_{j\geq 1}\frac{\|u' \varphi^{j-1}\|_{v_1}}{\|z^{j-1}\|_{v_{(2+\a)/p}}}<\infty~~~\mbox{and}~~~ \sup_{j\geq 1}\frac{\|u\varphi' \varphi^{j-1}\|_{v_1}}{\|z^{j-1}\|_{v_{(2+\a+p)/p}}}<\infty,
$$
respectively.   Since $I_uf(0)=0, J_uf(0)=0,$
  $$
\Big(  I_u(\varphi^j)(z)\Big)'=ju(z)\varphi'(z) \varphi^{j-1}(z),   \Big(  J_u(\varphi^{j-1})(z)\Big)'=u'(z) \varphi^{j-1}(z),
  $$
  by Lemma 3.3, we see that
  $uC_\varphi: A^p_\alpha  \to \mathcal{B}   $ is bounded if and only if
 \begr
\sup_{j\geq 1} j^{(2+\a)/p}\|J_u(\varphi^{j-1})\|_{\mathcal{B}}&=& \sup_{j\geq 1}
j^{(2+\a)/p}\|u' \varphi^{j-1}\|_{v_1} \nonumber\\
&\approx& \sup_{j\geq 1}\frac{j^{(2+\a)/p}\|u' \varphi^{j-1}\|_{v_1}}{j^{(2+\a)/p}\|z^{j-1}\|_{v_{(2+\a)/p}}}<\infty 
\endr
and
\begr
\sup_{j\geq 1} j^{(2+\a)/p}\|I_u(\varphi^j)\|_{\mathcal{B} }&= & \sup_{j\geq 1} j^{(2+\a+p)/p}\|u\varphi' \varphi^{j-1}\|_{v_1} \nonumber\\
&\approx& \sup_{j\geq 1}\frac{j^{(2+\a+p)/p}\|u\varphi' \varphi^{j-1}\|_{v_1}}{j^{(2+\a+p)/p}\|z^{j-1}\|_{v_{(2+\a+p)/p}}}<\infty.
\endr
 The proof is complete. \endproof

\noindent{\bf Theorem 3.2.}  Let $1\leq p <\infty$, $\a>-1$, $u \in H(\mathbb{D})$ and  $\varphi\in S(\DD)$ such that the operator $uC_\varphi: A^p_\alpha  \to \mathcal{B}   $ is bounded. Then
 \begr   \| uC_\varphi\|_{e, A^p_\alpha  \to \mathcal{B}  }  \approx   \max\Big\{\limsup_{j\rightarrow\infty} j^{(2+\a)/p} \|I_u(\varphi^j)\|_{\mathcal{B} },
 \limsup_{j\rightarrow\infty} j^{(2+\a)/p}\|J_u(\varphi^{j-1})\|_{\mathcal{B} }\Big\}.  \nonumber \endr

\proof   By Theorem A and Lemma 3.1, $uC_\varphi: A^p_\alpha  \to \mathcal{B}  $ is bounded if and only if  the weighted composition operator $u'C_\varphi: H^\infty_{v_{(2+\a)/p}}\rightarrow  H^\infty_{v_1}$ is bounded and $u\varphi'C_\varphi: H^\infty_{v_{(2+\a+p)/p}}\rightarrow  H^\infty_{v_1}$ is bounded. By Lemmas 3.2 and  3.3, we get
 \begr
 \|u'C_\varphi\|_{e, H^\infty_{v_{(2+\a)/p}}\rightarrow  H^\infty_{v_1}}&=& \limsup_{j\rightarrow\infty}\frac{\|u' \varphi^{j-1}\|_{v_1}}{\|z^{j-1}\|_{v_{(2+\a)/p}}}=   \limsup_{j\rightarrow\infty}\frac{j^{(2+\a)/p}\|u' \varphi^{j-1}\|_{v_1}}{j^{(2+\a)/p}\|z^{j-1}\|_{v_{(2+\a)/p}}}  \nonumber\\
 &\approx& \limsup_{j\rightarrow\infty}j^{(2+\a)/p}\|u' \varphi^{j-1}\|_{v_1}  \nonumber\\
   &=& \limsup_{j\rightarrow\infty} j^{(2+\a)/p}\|J_u(\varphi^{j-1})\|_{\mathcal{B}}
 \endr
and
 \begr
 \|u\varphi'C_\varphi\|_{e, H^\infty_{v_{(2+\a+p)/p}}\rightarrow  H^\infty_{v_1}}&=& \limsup_{j\rightarrow\infty}\frac{\|u\varphi' \varphi^{j-1}\|_{v_1}}{\|z^{j-1}\|_{v_{(2+\a+p)/p}}} \nonumber\\
 &\approx& \limsup_{j\rightarrow\infty}j^{(2+\a+p)/p}\|u\varphi' \varphi^{j-1}\|_{v_1} \nonumber\\
  &=& \limsup_{j\rightarrow\infty}j^{(2+\a)/p}\|I_u(\varphi^j)\|_{\mathcal{B}}.
 \endr

{\bf The upper estimate.} From the fact   $(uC_\varphi f)'(z)=u'(z)f(\varphi(z))+u(z)\varphi'(z)f'(\varphi(z))$, it is easy to see  that
 \begr
 \| uC_\varphi\|_{e, A^p_\alpha  \to \mathcal{B} }  \leq  \|u'C_\varphi\|_{e, H^\infty_{v_{(2+\a)/p}}\rightarrow  H^\infty_{v_1}}+\|u\varphi'C_\varphi\|_{e, H^\infty_{v_{(2+\a+p)/p}}\rightarrow  H^\infty_{v_1}}.
 \endr
Then, by (3.5), (3.6) and (3.7) we get
 \begr   \| uC_\varphi\|_{e, A^p_\alpha  \to \mathcal{B}  }
&\lesssim & \limsup_{j\rightarrow\infty}j^{(2+\a)/p}\|I_u(\varphi^j)\|_{\mathcal{B} }  +\limsup_{j\rightarrow\infty} j^{(2+\a)/p}\|J_u(\varphi^{j-1})\|_{\mathcal{B}} \nonumber \\
&\lesssim &  \max\Big\{\limsup_{j\rightarrow\infty} j^{(2+\a)/p} \|I_u(\varphi^j)\|_{\mathcal{B} },
 \limsup_{j\rightarrow\infty} j^{(2+\a)/p}\|J_u(\varphi^{j-1})\|_{\mathcal{B} }\Big\}.  \nonumber \endr

{\bf The lower estimate.} From Theorem 2.1 and Lemma 3.1,  we have
  \begr   \| uC_\varphi \|_{e, A^p_\alpha  \to \mathcal{B}  } &\gtrsim &   P
 =     \|u'C_\varphi\|_{e, H^\infty_{(2+\a)/p}\rightarrow  H^\infty_{v_1}} \approx   \limsup_{j\rightarrow\infty}  j^{(2+\a)/p} \|J_u(\varphi^{j-1})\|_{\mathcal{B} }    \nonumber \endr
 and
   \begr   \| uC_\varphi\|_{e, A^p_\alpha  \to \mathcal{B} } &\gtrsim &Q
 = \|u \varphi' C_\varphi\|_{e, H^\infty_{v_{(2+\a+p)/p}}\rightarrow  H^\infty_{v_1}}  \approx   \limsup_{j\rightarrow\infty}j^{(2+\a)/p}\|I_u(\varphi^j)\|_{\mathcal{B} } .   \nonumber \endr
Therefore,
  \begr   \| uC_\varphi\|_{e, A^p_\alpha  \to \mathcal{B}  } \gtrsim  \max\Big\{\limsup_{j\rightarrow\infty} j^{(2+\a)/p} \|I_u(\varphi^j)\|_{\mathcal{B} }, \limsup_{j\rightarrow\infty} j^{(2+\a)/p}\|J_u(\varphi^{j-1})\|_{\mathcal{B} }\Big\}.  \nonumber \endr
This completes the proof.
\endproof

From Theorem 3.2, we immediately get the following result.
\msk

\noindent{\bf Theorem 3.3.}    Let    $1\leq p <\infty$, $\a>-1$, $u \in H(\mathbb{D})$ and  $\varphi\in S(\DD)$ such that $uC_\varphi: A^p_\alpha  \to \mathcal{B}  $ is bounded. Then the operator $uC_\varphi: A^p_\alpha  \to \mathcal{B}  $ is compact if and only if
  \begr   \limsup_{j\rightarrow\infty} j^{(2+\a)/p} \|I_u(\varphi^j)\|_{\mathcal{B} }=0~~~~~\mbox{and}~~~~~~
 \limsup_{j\rightarrow\infty} j^{(2+\a)/p}\|J_u(\varphi^{j-1})\|_{\mathcal{B} } =0.  \nonumber \endr

We end this section with  a new characterization of boundedness, compactness and essential norm of the operator  $uC_\varphi: H^p  \to \mathcal{B}   $,  which
follows from Theorems 3.1, 3.2 and 3.3 by taking the limit as $\a$ decreases to $-1$.\msk

\noindent{\bf Theorem 3.4.} {\it   Let    $1 \leq p <\infty$, $u \in H(\mathbb{D})$ and  $\varphi\in S(\DD)$. Then the following statements hold.

(a) The operator $uC_\varphi: H^p  \to \mathcal{B}   $ is bounded if and only if
  \begr\sup_{j\geq 1} j^{1/p} \|I_u(\varphi^j)\|_{\mathcal{B} } <\infty~~~~~~~\mbox{and}~~~~~~
 \sup_{j\geq 1} j^{1/p}\|J_u(\varphi^{j-1})\|_{\mathcal{B} }   <\infty. \nonumber \endr

 (b) If the operator $uC_\varphi: H^p \to \mathcal{B}   $ is bounded. Then   $uC_\varphi: H^p  \to \mathcal{B}  $ is compact if and only if
  \begr   \limsup_{j\rightarrow\infty} j^{1/p} \|I_u(\varphi^j)\|_{\mathcal{B} }=0~~~~~\mbox{and}~~~~~~
 \limsup_{j\rightarrow\infty} j^{1/p}\|J_u(\varphi^{j-1})\|_{\mathcal{B} } =0.  \nonumber \endr
 Moreover
 \begr   \| uC_\varphi\|_{e, H^p  \to \mathcal{B}  }  \approx   \max\Big\{\limsup_{j\rightarrow\infty} j^{1/p} \|I_u(\varphi^j)\|_{\mathcal{B} },
 \limsup_{j\rightarrow\infty} j^{1/p}\|J_u(\varphi^{j-1})\|_{\mathcal{B} }\Big\}.  \nonumber \endr}

From the above results, we immediately get the following new characterization of the operator  $C_\varphi: A^p_\a ( \mbox{or} ~H^p) \to \mathcal{B}   $.\msk

\noindent{\bf Corollary 3.1.} {\it   Let $1 \leq p <\infty$, $\a>-1$  and  $\varphi\in S(\DD)$. Then the following statements hold.

(a) The operator $C_\varphi: A^p_\a  \to \mathcal{B}   $ is bounded if and only if
  $\sup_{j\geq 1} j^{(\a+2)/p} \| \varphi^j \|_{\mathcal{B} } <\infty . $

 (b) If the operator $C_\varphi: A^p_\a \to \mathcal{B}   $ is bounded. Then   $C_\varphi: A^p_\a \to \mathcal{B}  $ is compact if and only if
  $  \limsup_{j\rightarrow\infty} j^{(\a+2)/p} \| \varphi^j \|_{\mathcal{B} }=0.$
 Moreover,
 \begr   \| C_\varphi\|_{e, A^p_\a  \to \mathcal{B}  }  \approx   \limsup_{j\rightarrow\infty} j^{(\a+2)/p} \| \varphi^j \|_{\mathcal{B} } .  \nonumber \endr}

\noindent{\bf Corollary 3.2.} {\it   Let $1 \leq p <\infty$  and  $\varphi\in S(\DD)$. Then the following statements hold.

(a) The operator $C_\varphi: H^p  \to \mathcal{B}   $ is bounded if and only if
  $\sup_{j\geq 1} j^{1/p} \| \varphi^j \|_{\mathcal{B} } <\infty . $

 (b) If the operator $C_\varphi: H^p \to \mathcal{B}   $ is bounded, then   $C_\varphi: H^p  \to \mathcal{B}  $ is compact if and only if
$ \limsup_{j\rightarrow\infty} j^{1/p} \| \varphi^j \|_{\mathcal{B} }=0.  $
 Moreover
 \begr   \| C_\varphi\|_{e, H^p  \to \mathcal{B}  }  \approx   \limsup_{j\rightarrow\infty} j^{1/p} \| \varphi^j \|_{\mathcal{B} } .  \nonumber \endr}

{\bf Acknowledgement.}   This project was partially supported by the Macao Science and Technology Development Fund(No.098/2013/A3), NSF
of Guangdong Province(No.S2013010011978) and NNSF of China(No.11471143).

\end{document}